# On the posterior distribution of classes of random means

LANCELOT F. JAMES[1], ANTONIO LIJOI[2] and IGOR PRÜNSTER[3]

[1]*Department of Information Systems, Business Statistics and Operations Management, Hong Kong University of Science and Technology, Clear Water Bay, Kowloon, Hong Kong.*
*E-mail: lancelot@ust.hk*

[2]*Dipartimento Economia Politica e Metodi Quantitatavi, Università degli Studi di Pavia, via San Felice 5, 27100 Pavia, and CNR-IMATI, Milano, Italy. E-mail: lijoi@unipv.it*

[3]*Dipartimento di Statistica e Matematica Applicata, Collegio Carlo Alberto and ICER, Università degli Studi di Torino, Piazza Arbarello 8, 10122 Torino, Italy. E-mail: igor@econ.unito.it*

The study of properties of mean functionals of random probability measures is an important area of research in the theory of Bayesian nonparametric statistics. Many results are now known for random Dirichlet means, but little is known, especially in terms of posterior distributions, for classes of priors beyond the Dirichlet process. In this paper, we consider normalized random measures with independent increments (NRMI's) and mixtures of NRMI. In both cases, we are able to provide exact expressions for the posterior distribution of their means. These general results are then specialized, leading to distributional results for means of two important particular cases of NRMI's and also of the two-parameter Poisson–Dirichlet process.

*Keywords:* Bayesian nonparametrics; completely random measures; means of random probability measures; normalized random measures; Poisson–Dirichlet process; posterior distribution; species sampling models

## 1. Introduction

Bayesian nonparametrics has recently undergone major development which has led to the proposal of a variety of new classes of prior distributions, as well as allowing the concrete application of nonparametric models to problems in, for example, biology, medicine, economics and ecology. While there is a vast literature on computational issues related to Bayesian nonparametric procedures, there is a dearth of analytical results, mainly due to the difficulties of studying distributions on infinite-dimensional spaces. Indeed, given a nonparametric prior, a natural statistical object to analyze is the mean functional: for instance, in the context of survival analysis, the mean takes on the interpretation of (random) expected lifetime. However, such an analysis seemed to be prohibitive until the pioneering contributions of Cifarelli and Regazzini [4, 5] who set up a general theory for the study of Dirichlet process means and also derived the remarkable *Markov–Krein* or *Cifarelli–Regazzini* identity. Since then, much attention has been focused on means of the Dirichlet process. Among other contributions, we mention [7, 8, 11, 13, 16, 17, 19, 31, 39, 43]. Recently the first results concerning nonlinear functionals of the Dirichlet process, such as the variance functional, have appeared in the literature (see [6, 34, 39]). Another line of research has dealt with mean functionals for priors different from the Dirichlet process (see [10, 14, 22, 33, 40, 42]). The study of means also highlights the interplay with other areas of mathematics, such as special functions [31], excursion theory [22, 42] and mathematical physics [38].

While some results concerning the prior distribution of means for classes of priors more general than the Dirichlet process are known, no exact result is known for their posterior distribution. Indeed, in [40], normalized random measures with independent increments (NRMI's), whose construction is recalled in Definition 2.1, were considered: in addition to results for the existence and the exact prior distribution of their means, an approximation for their posterior density were achieved. These results were then extended in [33] to mixtures of NRMI, leading to an approximation of the posterior mean of mixtures of NRMI and to the exact expression for the special case of the mixture of Dirichlet process. These two papers [33, 40] represent the starting point of our work and we aim to develop and complete these results: indeed, we are able to provide exact expressions for the posterior distributions of both means of NRMI's and means of mixtures of NRMI.

The outline of the paper is as follows. In Section 2, the basic concepts are introduced and preliminary results recalled. In Section 3.1, we determine exact expressions for the posterior distribution of means of NRMI's, whereas in Section 3.2, general formulae for posterior means of mixtures of NRMI are obtained. Section 4 is devoted to the study of means of particular NRMI's of statistical relevance, namely the extended gamma NRMI and the generalized gamma NRMI, and the main result of [33] is recovered as a corollary. Moreover, our results for the generalized gamma NRMI are exploited to derive a new expression for the distribution of a mean of the two-parameter Poisson–Dirichlet process. Proofs are deferred to the Appendix.

## 2. Preliminaries and basic definitions

We first recall the concept of completely random measure, due to Kingman [26]. Let $(\mathbb{X}, \mathscr{X})$ be a Polish space endowed with the Borel $\sigma$-field and $(\mathbb{M}, \mathscr{B}(\mathbb{M}))$ be the space of boundedly finite measures on $\mathbb{X}$, with $\mathscr{B}(\mathbb{M})$ denoting the corresponding Borel $\sigma$-algebra. Let $\tilde{\mu}$ be a random element defined on $(\Omega, \mathscr{F}, \mathbb{P})$ and with values in $(\mathbb{M}, \mathscr{B}(\mathbb{M}))$ such that for any collection of disjoint sets in $\mathscr{X}$, $A_1, A_2, \ldots$, the random variables $\tilde{\mu}(A_1), \tilde{\mu}(A_2), \ldots$ are mutually independent. Then, $\tilde{\mu}$ is a *completely random measure* (CRM) on $(\mathbb{X}, \mathscr{X})$. A CRM can always be represented as a linear functional of a Poisson random measure. In particular, define $\mathscr{H}_\nu$ to be the space of measurable functions $h : \mathbb{X} \to \mathbb{R}^+$ such that $\int_{\mathbb{R}^+ \times \mathbb{X}} [1 - \mathrm{e}^{-vh(x)}] \nu(\mathrm{d}v, \mathrm{d}x) < \infty$, where $\nu$ stands for the intensity of the Poisson random measure underlying $\tilde{\mu}$, which must satisfy the integrability condition $\int_{\mathbb{R}^+} \min\{v, 1\} \nu(\mathrm{d}v, \mathrm{d}x) < +\infty$ for almost all $x \in \mathbb{X}$. Then, $\tilde{\mu}$ is uniquely characterized by its *Laplace functional* which, for any $h$ in $\mathscr{H}_\nu$, is given by

$$\mathbb{E}[\mathrm{e}^{-\int_\mathbb{X} h(x)\tilde{\mu}(\mathrm{d}x)}] = \mathrm{e}^{-\int_{\mathbb{R}^+\times\mathbb{X}}[1-\mathrm{e}^{-vh(x)}]\nu(\mathrm{d}v,\mathrm{d}x)}. \tag{1}$$

Throughout the paper, we define, for any real- or complex-valued function $g$ defined on $\mathbb{X}$, the functional $\psi(g) := \int_{\mathbb{R}^+\times\mathbb{X}}[1-\mathrm{e}^{-vg(x)}]\nu(\mathrm{d}v,\mathrm{d}x)$. Moreover, let $\mathbf{1}:\mathbb{X}\to\mathbb{R}$ be the function identically equal to 1, namely $\mathbf{1}(x)=1$ for any $x$ in $\mathbb{X}$. See [28] for an exhaustive account on CRM's. The representation in (1) establishes that $\tilde{\mu}$ is characterized by the corresponding intensity $\nu$. Letting $\alpha$ be a $\sigma$-finite measure on $\mathbb{X}$, we can always write

$$\nu(\mathrm{d}v,\mathrm{d}x) = \rho(\mathrm{d}v|x)\alpha(\mathrm{d}x), \tag{2}$$

where $\rho$ is a measurable kernel such that $\rho(\cdot|x)$ is a $\sigma$-finite measure on $\mathscr{B}(\mathbb{R}^+)$ for any $x$ in $\mathbb{X}$. Such a disintegration is ensured by [25], Theorem 15.3.3. If (2) simplifies to

$$\nu(\mathrm{d}v,\mathrm{d}x) = \rho(\mathrm{d}v)\alpha(\mathrm{d}x), \tag{3}$$

then we say that $\tilde{\mu}$ is a *homogeneous* CRM, whereas if this is not the case, then $\tilde{\mu}$ will be termed a *non-homogeneous* CRM. In the following, we will assume $\alpha$ is non-atomic.

Since the aim is to define random probability measures by means of normalization of completely random measures, the total mass $\tilde{\mu}(\mathbb{X})$ needs to be finite and positive, almost surely. As shown in [40], this happens if $\nu(\mathbb{R}^+\times\mathbb{X}) = +\infty$ and the Laplace exponent $\psi(\lambda) < +\infty$ for any $\lambda > 0$, respectively.

**Definition 2.1.** *Let $\tilde{\mu}$ be a CRM on $(\mathbb{X},\mathscr{X})$ such that $\nu(\mathbb{R}^+\times\mathbb{X}) = +\infty$ and $\psi(\lambda) < +\infty$ for any $\lambda > 0$. The random probability measure*

$$\tilde{P}(\cdot) = \frac{\tilde{\mu}(\cdot)}{\tilde{\mu}(\mathbb{X})} \tag{4}$$

*is then termed an NRMI.*

Strictly speaking, the random probability measure in (4) is a normalized CRM and reduces to an NRMI when $\mathbb{X} = \mathbb{R}$. Nonetheless, we prefer to keep the acronym NRMI introduced in [40]. According to the decomposition of the intensity $\nu$ described in (2) and (3), we will distinguish between non-homogeneous and homogeneous NRMI's. Several priors used in Bayesian nonparametric inference can be defined as in (4). For instance, as already noted by [12], the Dirichlet process can be recovered as an NRMI based on the gamma CRM for which $\nu(\mathrm{d}v,\mathrm{d}x) = v^{-1}\mathrm{e}^{-v}\mathrm{d}v\,\alpha(\mathrm{d}x)$. Other examples include the normalized stable process [27] and the normalized inverse Gaussian process [29].

Nowadays, the most common use of Bayesian nonparametric procedures is within hierarchical mixtures: letting $\mathbb{Y}$ be a Polish space equipped with the Borel $\sigma$-algebra $\mathscr{Y}$, one defines a random density (absolutely continuous with respect to some $\sigma$-finite measure $\lambda$ on $\mathbb{Y}$) driven by a random discrete distribution, that is,

$$\tilde{f}(y) = \int_\mathbb{X} k(y,x)\tilde{P}(\mathrm{d}x), \tag{5}$$

where $k$ is a density function on $\mathbb{Y}$ indexed by some parameter with values in $\mathbb{X}$. A typical choice for $k$ is represented by the density function of the normal distribution: in such a case, $\tilde{P}$ controls the means (and possibly also the variances) of the random mixture density. This approach is due to Lo [32], who defined a random density as in (5) with $\tilde{P}$ being the Dirichlet process; this model is now commonly referred to as mixture of Dirichlet process (MDP). Recently, various contributions have focused on replacing the Dirichlet process in (5) with alternative random probability measures, which yield interesting behaviors, especially in terms of the induced clustering mechanism; see, for example, [18, 20, 29].

The present paper is focused on linear functionals of NRMI's, namely on random quantities of the type $\tilde{P}(g) := \int_{\mathbb{X}} g(x)\tilde{P}(\mathrm{d}x)$, where $g:\mathbb{X} \to \mathbb{R}$ is any measurable function such that

$$\psi(t|g|) = \int_{\mathbb{X}\times\mathbb{R}^+} (1 - \mathrm{e}^{-tv|g(x)|})\rho(\mathrm{d}v|x)\alpha(\mathrm{d}x) < +\infty \qquad \forall t > 0. \tag{6}$$

By [40], Proposition 1, condition (6) is necessary and sufficient for $\tilde{P}(|g|)$ to be a.s. finite. In the sequel, we always assume that (6) holds true. An exact analytic expression for $\mathbb{F}$, the distribution function of $\tilde{P}(g)$, is given in [40], Proposition 2. We will also examine means of a mixture of NRMI, that is

$$\tilde{Q}(g) := \int_{\mathbb{Y}} g(y)\tilde{f}(y)\lambda(\mathrm{d}y) = \int_{\mathbb{X}} h(x)\tilde{P}(\mathrm{d}x) = \tilde{P}(h), \tag{7}$$

where $\tilde{Q}$ stands for the random probability measure associated with $\tilde{f}$, defined as in (5), and $h(x) = \int_{\mathbb{Y}} g(y)k(y,x)\lambda(\mathrm{d}y)$. Hence, as shown in [33], the necessary and sufficient condition for $\tilde{Q}(g)$ being a.s. finite becomes

$$\psi(th^*) = \int_{\mathbb{X}\times\mathbb{R}^+} (1 - \mathrm{e}^{-tvh^*(x)})\rho(\mathrm{d}v|x)\alpha(\mathrm{d}x) < +\infty \qquad \forall t > 0, \tag{8}$$

with $h^*(x) = \int_{\mathbb{Y}} |g(y)|k(y,x)\lambda(\mathrm{d}y)$. The evaluation of the prior distribution of the mean then follows in a straightforward way [33], Proposition 2. In the following, when we consider means of mixtures of NRMI, as in (7), we will tacitly suppose that $g$ verifies condition (8).

## 3. Posterior distribution of means

### 3.1. Means of NRMI's

We first focus attention on posterior distributions of means of NRMI's. Let $(X_n)_{n\geq 1}$ be a sequence of exchangeable observations, defined on $(\Omega, \mathscr{F}, \mathbb{P})$ and with values in $\mathbb{X}$, such that, given an NRMI $\tilde{P}$, the $X_i$'s are i.i.d. with distribution $\tilde{P}$, that is, for any $B_i \in \mathscr{X}$,

$i = 1, \ldots, n$ and $n \geq 1$,

$$\mathbb{P}[X_1 \in B_1, \ldots, X_n \in B_n | \tilde{P}] = \prod_{i=1}^{n} \tilde{P}(B_i). \tag{9}$$

Moreover, let $\mathbf{X} = (X_1, \ldots, X_n)$. It is clear that one can always represent $\mathbf{X}$ as $(\mathbf{X}^*, \boldsymbol{\pi})$, where $\mathbf{X}^* = (X_1^*, \ldots, X_{n(\boldsymbol{\pi})}^*)$ denotes the distinct observations within the sample and $\boldsymbol{\pi} = \{C_1, \ldots, C_{n(\boldsymbol{\pi})}\}$ stands for the corresponding partition of the integers $\{1, \ldots, n\}$ recording which observations within the sample are equal, that is, $C_j = \{i : X_i = X_j^*\}$. The number of elements in the $j$th set of the partition is indicated by $n_j$, for $j = 1, \ldots, n(\boldsymbol{\pi})$, so that $\sum_{j=1}^{n(\boldsymbol{\pi})} n_j = n$.

At this point, it is useful to recall the posterior characterization of NRMI's given in [23]. For any pair of random elements $Z$ and $W$ defined on $(\Omega, \mathscr{F}, \mathbb{P})$, we use the symbol $Z^{(W)}$ to denote a random element on $(\Omega, \mathscr{F}, \mathbb{P})$ whose distribution coincides with a regular conditional distribution of $Z$, given $W$. Now, introduce a latent variable, denoted by $U_n$, whose conditional distribution, given $\mathbf{X}$, admits a density function (with respect to the Lebesgue measure on $\mathbb{R}$) coinciding with

$$f_{U_n}^{\mathbf{X}}(u) \propto u^{n-1} \prod_{i=1}^{n(\boldsymbol{\pi})} \tau_{n_i}(u | X_i^*) e^{-\psi(u)}, \tag{10}$$

where

$$\tau_{n_i}(u | X_i^*) = \int_{\mathbb{R}^+} s^{n_i} e^{-us} \rho(\mathrm{d}s | X_i^*) \tag{11}$$

for $i = 1, \ldots, n(\boldsymbol{\pi})$. Indeed, the posterior distribution, given $\mathbf{X}$, of the CRM $\tilde{\mu}$ defining an NRMI (4) is a mixture with respect to the distribution of the latent variable $U_n$. Specifically, $\tilde{\mu}^{(U_n, \mathbf{X})} \stackrel{d}{=} \tilde{\mu}^{(U_n)} + \sum_{i=1}^{n(\boldsymbol{\pi})} J_i^{(U_n, \mathbf{X})} \delta_{X_i^*}$, where: $\tilde{\mu}^{(U_n)}$ is a CRM with intensity

$$\nu^{(U_n)}(\mathrm{d}s, \mathrm{d}x) = e^{-U_n s} \rho(\mathrm{d}s | x) \alpha(\mathrm{d}x); \tag{12}$$

the $X_i^*$'s are the fixed points of discontinuity; the $J_i^{(U_n, \mathbf{X})}$'s are the corresponding jumps, which are mutually independent and independent from $\tilde{\mu}^{(U_n)}$, and whose density is given by

$$f_{J_i}^{(U_n, \mathbf{X})}(s) \propto s^{n_i} e^{-U_n s} \rho(\mathrm{d}s | X_i^*). \tag{13}$$

See [23], Theorem 1, for details.

We are now in a position to provide the exact posterior distribution of means of an NRMI. Note that the results hold for any function $g$ and any NRMI (identified by means of its Poisson intensity (2)) which lead to an a.s. finite mean (6). In what follows, we agree to denote by $\mathrm{Im}(z)$ and $\mathrm{Re}(z)$ the imaginary and real parts, respectively, of the complex number $z$. Moreover, $\psi^{(u)}$ and $J_r^{(u, \mathbf{X})}$ are the Laplace exponent of the CRM defined by (12) and the jumps whose density is given by (13), respectively, with $U_n = u$.

**Theorem 3.1.** *Let $\tilde{P}$ be an NRMI. The posterior distribution of $\tilde{P}(g)$, given $\mathbf{X}$, is then absolutely continuous with respect to the Lebesgue measure on $\mathbb{R}$ and a posterior density function is given by*

$$\rho^{\mathbf{X}}(\sigma; g) = \begin{cases} \int_0^\infty \mathrm{Re}\{\chi_g(t,\sigma)\}\,\mathrm{d}t, & \text{if } n=1, \\ (-1)^{p+1} \int_{-\infty}^\sigma \int_0^\infty [(\sigma-z)t]^{n-1} \mathrm{Im}\{\chi_g(t,z)\}\,\mathrm{d}t\,\mathrm{d}z, & \text{if } n=2p, \\ (-1)^p \int_{-\infty}^\sigma \int_0^\infty [(\sigma-z)t]^{n-1} \mathrm{Re}\{\chi_g(t,z)\}\,\mathrm{d}t\,\mathrm{d}z, & \text{if } n=2p+1, \end{cases} \quad (14)$$

*where $p \geq 1$,*

$$\chi_g(t,z) = \frac{\mathrm{e}^{-\psi(-\mathrm{i}t(g-z\mathbf{1}))} \prod_{j=1}^{n(\boldsymbol{\pi})} \kappa_{n_j}(\mathrm{i}t[g(X_j^*)-z]|X_j^*)}{\pi \int_0^{+\infty} u^{n-1} [\prod_{j=1}^{n(\boldsymbol{\pi})} \tau_{n_j}(u|X_j^*)] \mathrm{e}^{-\psi(u\mathbf{1})}\,\mathrm{d}u}, \quad (15)$$

$\kappa_{n_j}(\mathrm{i}t[g(X_j^*)-z]|X_j^*) = \int_0^{+\infty} v^{n_j} \mathrm{e}^{\mathrm{i}tv(g(X_j^*)-z)} \rho(\mathrm{d}v|X_j^*)$ and $\tau_{n_j}(u|X_j^*)$ *is as in (11), for $j=1,\ldots,n(\boldsymbol{\pi})$. Moreover, the posterior cumulative distribution function of $\tilde{P}(g)$, given $\mathbf{X}$, admits the representation*

$$\mathbb{F}^{\mathbf{X}}(\sigma; g) = \frac{1}{2} - \frac{1}{\pi} \lim_{T \to +\infty} \int_0^T \frac{1}{t} \int_0^{+\infty} \zeta_g(\sigma; u,t) f_{U_n}^{\mathbf{X}}(u)\,\mathrm{d}u\,\mathrm{d}t, \quad (16)$$

*where*

$$\zeta_g(\sigma; u,t) := \mathrm{Im}\{\mathrm{e}^{-\psi^{(u)}(-\mathrm{i}t(g-\sigma\mathbf{1}))} \mathbb{E}[\mathrm{e}^{\mathrm{i}t \sum_{r=1}^{n(\boldsymbol{\pi})}(g(X_r^*)-\sigma)J_r^{(u,\mathbf{X})}}]\} \quad (17)$$

*and $f_{U_n}^{\mathbf{X}}$ is the density of the latent variable $U_n$ given in (10).*

### 3.2. Means of mixtures of NRMI

Before dealing with the posterior distribution of means of mixtures of NRMI, it is worth recalling that a popular way of representing (5) is as a hierarchical mixture:

$$\begin{aligned} Y_i | X_i &\overset{\text{i.n.d.}}{\sim} k(\cdot, X_i), \quad i=1,\ldots,n, \\ X_i | \tilde{P} &\overset{\text{i.i.d.}}{\sim} \tilde{P}, \\ \tilde{P} &\sim \mathscr{P}, \end{aligned} \quad (18)$$

where $(X_i)_{i \geq 1}$ is a sequence of latent variables with values in $\mathbb{X}$ and $\mathscr{P}$ is the law of the NRMI $\tilde{P}$. The notation adopted in the description of model (18) is standard in Bayesian statistics and means that, given $\tilde{P}$, the observations $Y_i$ are independent and identically distributed with random density $\tilde{f}(y) = \int_{\mathbb{X}} k(y,x) \tilde{P}(\mathrm{d}x)$. Now, since the background driving NRMI is a.s. discrete, it will induce ties among the latent variables, which are denoted, as previously, by $X_1^*, \ldots, X_{n(\boldsymbol{\pi})}^*$ and $\boldsymbol{\pi} = \{C_1, \ldots, C_{n(\boldsymbol{\pi})}\}$ indicates the corresponding partition of the integers $\{1, \ldots, n\}$ with $C_j = \{i : X_i = X_j^*\}$. This notation allows us, given a partition, to keep a record of which latent variable each of the observations $Y_i$, for $i = 1, \ldots, n$, is assigned to.

Here, we aim to derive the posterior distribution of $\tilde{Q}(g)$, given $\mathbf{Y}$. By virtue of (7), this can also be seen as the distribution of $\int_{\mathbb{X}} h(x) \tilde{P}(\mathrm{d}x)$, given $\mathbf{Y}$. The next theorem yields an exact expression for a density of such a posterior distribution. It is worth noting that the results hold for any function $g$, $k$ and NRMI which lead to a finite mean (8): the expressions are given in terms of $h$ (which is in turn, defined as in (7)), $k$ and the Poisson intensity (2) corresponding to any NRMI.

**Theorem 3.2.** *Let $\tilde{f}$ be a mixture of NRMI, as in (5). The posterior distribution of $\tilde{Q}(g)$, given $\mathbf{Y}$, is then absolutely continuous with respect to the Lebesgue measure on $\mathbb{R}$ and the corresponding posterior density function is given by*

$$\phi^{\mathbf{Y}}(\sigma; g) = \begin{cases} \int_0^\infty \mathrm{Re}\{\xi_h(t, \sigma)\} \, \mathrm{d}t, & \text{if } n = 1, \\ (-1)^{p+1} \int_{-\infty}^{\sigma} \int_0^\infty [(\sigma - z)t]^{n-1} \mathrm{Im}\{\xi_h(t, z)\} \, \mathrm{d}t \, \mathrm{d}z, & \text{if } n = 2p, \\ (-1)^p \int_{-\infty}^{\sigma} \int_{-\infty}^{\sigma} \int_0^\infty [(\sigma - z)t]^{n-1} \mathrm{Re}\{\xi_h(t, z)\} \, \mathrm{d}t \, \mathrm{d}z, & \text{if } n = 2p+1, \end{cases}$$

(19)

*where*

$$\xi_h(t, z) = \frac{\mathrm{e}^{-\psi(-\mathrm{i}t(h - z\mathbf{1}))} [\sum_{\boldsymbol{\pi}} \prod_{j=1}^{n(\boldsymbol{\pi})} \int_{\mathbb{X}} \kappa_{n_j}(\mathrm{i}t[h(x) - z] | x) \prod_{i \in C_j} k(Y_i, x) \alpha(\mathrm{d}x)]}{\pi \int_{\mathbb{R}^+} u^{n-1} \mathrm{e}^{-\psi(u\mathbf{1})} [\sum_{\boldsymbol{\pi}} \prod_{j=1}^{n(\boldsymbol{\pi})} \int_{\mathbb{X}} \prod_{i \in C_j} k(Y_i, x) \tau_{n_j}(u|x) \alpha(\mathrm{d}x)] \, \mathrm{d}u}, \quad (20)$$

*having set $h(x) = \int_{\mathbb{Y}} g(y) k(y, x) \lambda(\mathrm{d}y)$ and $\kappa_{n_j}(\mathrm{i}t[h(x) - z]|x) = \int_0^{+\infty} v^{n_j} \mathrm{e}^{\mathrm{i}tv(h(x)-z)} \rho(\mathrm{d}v|x)$ for $j = 1, \ldots, n(\boldsymbol{\pi})$. Moreover, the posterior cumulative distribution function of $\tilde{Q}(g)$, given $\mathbf{Y}$, can be represented as*

$$\mathbb{G}^{\mathbf{Y}}(\sigma; g) = \frac{1}{2} - \frac{1}{\pi} \lim_{T \to +\infty} \int_0^T \frac{1}{t} \mathrm{Im}\{\zeta_h(\sigma; t)\} \, \mathrm{d}t, \quad (21)$$

*where $p \geq 1$,*

$$\zeta_h(\sigma; t) = \int_{\mathbb{R}^+} u^{n-1} \mathrm{e}^{-\psi(-\mathrm{i}t(g - \sigma\mathbf{1}) + u\mathbf{1})}$$

$$\times \left[ \sum_{\boldsymbol{\pi}} \prod_{j=1}^{n(\boldsymbol{\pi})} \int_{\mathbb{X}} \mathbb{E}(\mathrm{e}^{\mathrm{i}t(g(x) - \sigma\mathbf{1}) J_j^{(u, \mathbf{X})}}) \prod_{i \in C_j} k(Y_i, x) \tau_{n_j}(u|x) \alpha(\mathrm{d}x) \right] \mathrm{d}u \quad (22)$$

$$\times \left( \int_{\mathbb{R}^+} u^{n-1} \mathrm{e}^{-\psi(u\mathbf{1})} \left[ \sum_{\boldsymbol{\pi}} \prod_{j=1}^{n(\boldsymbol{\pi})} \int_{\mathbb{X}} \prod_{i \in C_j} k(Y_i, x) \tau_{n_j}(u|x) \alpha(\mathrm{d}x) \right] \mathrm{d}u \right)^{-1},$$

where the jumps $J_j^{(U_n,\mathbf{X})}$'s have density given by (13).

## 4. Applications

Let us now consider, in detail, two specific cases of statistical relevance, involving two important CRM's, namely, generalized gamma CRM [2] and extended gamma CRM [9]. Both have been found to have many applications in survival analysis, spatial statistics, mixture models and spatio-temporal models.

### 4.1. Extended gamma NRMI

Here, we consider NRMI's derived from extended CRM's, which are characterized by the intensity

$$\nu(\mathrm{d}v, \mathrm{d}x) = \frac{\mathrm{e}^{-\beta(x)v}}{v}\,\mathrm{d}v\,\alpha(\mathrm{d}x),$$

where $\beta$ is a positive real-valued function. The corresponding NRMI in (4), to be termed an *extended gamma NRMI* with parameters $(\alpha(\cdot), \beta(\cdot))$, is well defined if $\alpha$ and $\beta$ are such that $\int_{\mathbb{X}} \log(1 + t[\beta(x)]^{-1})\alpha(\mathrm{d}x) < +\infty$ for every $t > 0$.

As for distributional properties of means of extended gamma NRMI's, from (6), it follows that $\tilde{P}(g)$ is finite if and only if $\int_{\mathbb{X}} \log(1 + tg(x)[\beta(x)]^{-1})\alpha(\mathrm{d}x) < +\infty$ for every $t > 0$, which coincides, except for the factor of non-homogeneity $\beta$, with the well-known condition given in [11] for the Dirichlet process. Moreover, the prior distribution of a mean is given by

$$\mathbb{F}(\sigma;g) = \frac{1}{2} - \frac{1}{\pi}\int_0^{+\infty} \frac{1}{t}\mathrm{Im}(\mathrm{e}^{-\int_{\mathbb{X}} \log(1-\mathrm{i}t[\beta(x)]^{-1}(g(x)-\sigma))\alpha(\mathrm{d}x)})\,\mathrm{d}t,$$

having applied Proposition 2 in [40] and shown that $\mathbb{F}$ is absolutely continuous with respect to the Lebesgue measure on $\mathbb{R}$. We now provide expressions for the posterior density function and posterior distribution function of a mean of an extended gamma NRMI.

**Proposition 4.1.** *Let $\tilde{P}$ be an extended gamma NRMI. The posterior density function of $\tilde{P}(g)$, given $\mathbf{X}$, is then of the form (14) with*

$$\chi_g(t,z) = \frac{\mathrm{e}^{-\int_{\mathbb{X}} \log(\beta(x) - \mathrm{i}t(g(x)-z))\alpha_n^{\mathbf{X}}(\mathrm{d}x)}}{\pi \int_0^{+\infty} u^{n-1}\mathrm{e}^{-\int_{\mathbb{X}} \log(\beta(x)+u)\alpha_n^{\mathbf{X}}(\mathrm{d}x)}\,\mathrm{d}u},$$

*where $\alpha_n^{\mathbf{X}}(\cdot) := \alpha(\cdot) + \sum_{i=1}^{n(\boldsymbol{\pi})} n_i \delta_{X_i^*}(\cdot)$. Moreover, the posterior cumulative distribution function is given by*

$$\mathbb{F}^{\mathbf{X}}(\sigma;g) = \frac{1}{2} - \frac{\int_0^\infty \int_0^\infty t^{-1}u^{n-1}\mathrm{Im}\{\mathrm{e}^{-\int_{\mathbb{X}} \log(\beta(x)+u-\mathrm{i}t(g(x)-\sigma))\alpha_n^{\mathbf{X}}(\mathrm{d}x)}\}\,\mathrm{d}u\,\mathrm{d}t}{\pi \int_{\mathbb{R}^+} u^{n-1}\mathrm{e}^{-\int_{\mathbb{X}} \log(\beta(x)+u)\alpha_n^{\mathbf{X}}(\mathrm{d}x)}\,\mathrm{d}u}.$$

It is worth noting that the expressions in Proposition 4.1 are surprisingly simple, given the fact that they are exact distributions of functionals of a non-conjugate prior. Indeed, in terms of complexity, they are no more involved than the known expressions for Dirichlet means [4, 5, 39]. For illustrative purposes, suppose $n = 1$ and let $\beta(x) = \beta_1 \mathbb{1}_A(x) + \beta_2 \mathbb{1}_{A^c}(x)$, where $\beta_1 > \beta_2 > 0$ and $A \in \mathscr{X}$. In this case, it can be easily seen that the normalizing constant in $\chi_g(t,z)$ coincides with

$$\int_0^{+\infty} \frac{1}{[\beta_1+u]^{\alpha_1(A)}[\beta_2+u]^{\alpha_1(A^c)}}\,\mathrm{d}u = \frac{\beta_1^{-\alpha_1(A)}\beta_2^{1-\alpha_1(A^c)}}{\alpha(\mathbb{X})}\,{}_2F_1\!\left(\alpha_1(A),1;\alpha(\mathbb{X})+1;1-\frac{\beta_2}{\beta_1}\right),$$

where, for simplicity, we have set $\alpha_1 := \alpha_1^{\mathbf{X}}$ and ${}_2F_1$ is the Gauss hypergeometric function. Hence, by resorting to (14), one has that a posterior density of $\tilde{P}(g)$, given $X_1$, evaluated at $z \in [0,1]$ coincides with

$$\frac{\alpha(\mathbb{X})}{\pi} \frac{\beta_1^{\alpha_1(A)}\beta_2^{\alpha_1(A^c)-1}}{{}_2F_1(\alpha_1(A),1;\alpha(\mathbb{X})+1;1-\beta_2/\beta_1)}$$
$$\times \int_0^{+\infty} \exp\!\left\{-\frac{1}{2}\int_{\mathbb{X}} \log[\beta^2(x)+t^2(g(x)-z)^2]\alpha_1(\mathrm{d}x)\right\} \qquad (23)$$
$$\times \cos\!\left(\int_{\mathbb{X}} \arctan\frac{t(g(x)-z)}{\beta(x)}\alpha_1(\mathrm{d}x)\right)\mathrm{d}t.$$

It is worth emphasizing that, since this is a density function, one obtains an integral representation for the hypergeometric function ${}_2F_1(\alpha_1(A),1;\alpha(\mathbb{X})+1;1-\beta_2\beta_1^{-1})$ in terms of $g$, that is,

$${}_2F_1\!\left(\alpha_1(A),1;\alpha(\mathbb{X})+1;1-\frac{\beta_2}{\beta_1}\right)$$
$$= \frac{\alpha(\mathbb{X})}{\pi}\beta_1^{\alpha_1(A)}\beta_2^{\alpha_1(A^c)-1}\int_0^1\int_0^{+\infty} e^{-1/2\int_{\mathbb{X}}\log[\beta^2(x)+t^2(g(x)-z)^2]\alpha_1(\mathrm{d}x)}$$
$$\times \cos\!\left(\int_{\mathbb{X}}\arctan\frac{t(g(x)-z)}{\beta(x)}\alpha_1(\mathrm{d}x)\right)\mathrm{d}t\,\mathrm{d}z.$$

If we further suppose that $g = \mathbb{1}_A$, then (23) provides a posterior density for $\tilde{P}(A)$, given $X_1$. In order to obtain a simplified expression, suppose that $A$ and $\alpha$ are such that $\alpha(A) = \alpha(A^c) = 1$ and $X_1 \in A$ so that $\alpha_1(A) = 2 = 1 + \alpha_1(A^c)$. The corresponding posterior density for $\tilde{P}(A)$ is then of the form

$$\frac{2\beta_1^2}{{}_2F_1(2,1;3;1-\beta_2/\beta_1)}\frac{z}{[\beta_1 z+\beta_2(1-z)]^2}\mathbb{1}_{[0,1]}(z). \qquad (24)$$

Details of the determination of (24) are given in the Appendix. Also, note that ${}_2F_1(2,1;3;1-\beta_2/\beta_1) = 2\beta_1^2\int_0^1 z[\beta_1 z+\beta_2(1-z)]^{-2}\,\mathrm{d}z = 2\beta_1^2\log(\beta_1/\beta_2)(\beta_1-\beta_2)^{-2} -$

$2\beta_1(\beta_1 - \beta_2)^{-1}$. Since the (prior) probability distribution function of $\tilde{P}(A)$ is

$$\mathbb{F}(\sigma; \mathbb{1}_A) = \frac{\beta_1 \sigma}{\beta_1 \sigma + \beta_2(1-\sigma)} \mathbb{1}_{[0,1]}(\sigma) + \mathbb{1}_{(1,\infty)}(\sigma)$$

and the corresponding density is $\beta_1 \beta_2 [\beta_1 \sigma + \beta_2(1-\sigma)]^{-2} \mathbb{1}_{[0,1]}(\sigma)$, one finds out that having observed $X_1 \in A$, the probability mass of $\tilde{P}(A)$ in a neighborhood of 0 decreases, while it increases in a neighborhood of 1. The opposite phenomenon is observed when $X_1$ is not in $A$.

We now move on to considering posterior distributions of means of mixtures of the extended Gamma NRMI. Recall that, in such a case, the data $\mathbf{Y}$ are i.i.d., given $\tilde{Q}$, where $\tilde{Q}$ is the random probability measure corresponding to the mixture density (5) driven by an extended gamma NRMI. As shown in [33], a mean functional $\tilde{Q}(g)$ is finite if and only if for every $t > 0$, one has $\int_{\mathbb{X}} \log(1 + th^*(x)[\beta(x)]^{-1}) \alpha(\mathrm{d}x) < +\infty$, where $h^*(x) = \int_{\mathbb{Y}} |g(y)| k(y,x) \lambda(\mathrm{d}y)$. The prior distribution function is then of the form

$$\mathbb{F}(\sigma; g) = \frac{1}{2} - \frac{1}{\pi} \int_0^{+\infty} \frac{1}{t} \mathrm{Im}(\mathrm{e}^{-\int_{\mathbb{X}} \log(1 - \mathrm{i}t(h(x)-\sigma)[\beta(x)]^{-1})\alpha(\mathrm{d}x)}) \, \mathrm{d}t,$$

with $h(x) = \int_{\mathbb{Y}} |g(y)| k(y,x) \lambda(\mathrm{d}y)$. The next proposition provides the corresponding posterior density and cumulative distribution function. For both, we provide two expressions: the first is in line with Theorem 3.1, whereas the second exploits a sort of quasi-conjugacy peculiar to gamma-like models (see also Remark 4.1 in Section 4.2).

**Proposition 4.2.** *Let $\tilde{Q}$ be the random probability measure associated with the random density (5) driven by an extended gamma NRMI. The posterior density function of a mean $\tilde{Q}(g)$, given $\mathbf{Y}$, is then of the form (19), with*

$$\xi_h(t, z)$$
$$= \mathrm{e}^{-\int_{\mathbb{X}} \log(\beta(s) - \mathrm{i}t(h(s)-z))\alpha(\mathrm{d}s)}$$
$$\times \sum_{\boldsymbol{\pi}} \prod_{j=1}^{n(\boldsymbol{\pi})} (n_j - 1)! \int_{\mathbb{X}} \mathrm{e}^{-\log(\beta(x) - \mathrm{i}t(h(x)-z))n_j} \prod_{i \in C_j} k(Y_i, x) \alpha(\mathrm{d}x) \qquad (25)$$
$$\times \left( \pi \int_{\mathbb{R}^+} u^{n-1} \mathrm{e}^{-\int_{\mathbb{X}} \log(\beta(s)+u)\alpha(\mathrm{d}s)} \right.$$
$$\left. \times \sum_{\boldsymbol{\pi}} \prod_{j=1}^{n(\boldsymbol{\pi})} (n_j - 1)! \int_{\mathbb{X}} \mathrm{e}^{-\log(\beta(x)+u)n_j} \prod_{i \in C_j} k(Y_i, x) \alpha(\mathrm{d}x) \, \mathrm{d}u \right)^{-1},$$

*or, alternatively, as (19), with*

$$\xi_h(t, z) \tag{26}$$

$$= \frac{\sum_{\boldsymbol{\pi}} \prod_{j=1}^{n(\boldsymbol{\pi})}(n_j-1)! \int_{\mathbb{X}^{n(\boldsymbol{\pi})}} \mathrm{e}^{-\int_{\mathbb{X}} \log(\beta(s)-\mathrm{i}t(h(s)-z))\alpha_n^{\mathbf{X}}(\mathrm{d}s)} \prod_{i\in C_j} k(Y_i,x_j)\alpha(\mathrm{d}x_j)}{\pi \sum_{\boldsymbol{\pi}} \prod_{j=1}^{n(\boldsymbol{\pi})}(n_j-1)! \int_{\mathbb{X}^{n(\boldsymbol{\pi})}} \int_{\mathbb{R}^+} u^{n-1} \mathrm{e}^{-\int_{\mathbb{X}} \log(\beta(s)+u)\alpha_n^{\mathbf{X}}}\,\mathrm{d}u \prod_{i\in C_j} k(Y_i,x_j)\alpha(\mathrm{d}x_j)},$$

where $\alpha_n^{\mathbf{X}}(\cdot) = \alpha(\cdot) + \sum_{j=1}^{n(\boldsymbol{\pi})} n_j \delta_{x_j}$ and, as previously, $h(x) = \int_{\mathbb{Y}} g(y)k(y,x)\lambda(\mathrm{d}y)$. Moreover, the posterior distribution function can be represented as (21), with

$$\zeta_h(\sigma;t)$$
$$= \int_{\mathbb{R}^+} u^{n-1} \mathrm{e}^{-\int_{\mathbb{X}} \log(u+\beta(s)-\mathrm{i}t(h(s)-\sigma))\alpha(\mathrm{d}s)} \tag{27}$$
$$\times \sum_{\boldsymbol{\pi}} \prod_{j=1}^{n(\boldsymbol{\pi})}(n_j-1)! \int_{\mathbb{X}} \mathrm{e}^{-\log(u+\beta(x)-\mathrm{i}t(h(x)-\sigma))n_j} \prod_{i\in C_j} k(Y_i,x)\alpha(\mathrm{d}x)\,\mathrm{d}u$$
$$\times \left(\pi \int_{\mathbb{R}^+} u^{n-1} \mathrm{e}^{-\int_{\mathbb{X}} \log(\beta(s)+u)\alpha(\mathrm{d}s)} \right.$$
$$\left. \times \sum_{\boldsymbol{\pi}} \prod_{j=1}^{n(\boldsymbol{\pi})}(n_j-1)! \int_{\mathbb{X}} \mathrm{e}^{-\log(\beta(x)+u)n_j} \prod_{i\in C_j} k(Y_i,x)\alpha(\mathrm{d}x)\,\mathrm{d}u \right)^{-1},$$

or, alternatively, as (21) with

$$\zeta_h(\sigma;t)$$
$$= \sum_{\boldsymbol{\pi}} \prod_{j=1}^{n(\boldsymbol{\pi})}(n_j-1)! \int_{\mathbb{X}^{n(\boldsymbol{\pi})}} \int_{\mathbb{R}^+} u^{n-1} \mathrm{e}^{-\int_{\mathbb{X}} \log(u+\beta(s)-\mathrm{i}t(h(s)-\sigma))\alpha_n^{\mathbf{X}}}\,\mathrm{d}u$$
$$\times \prod_{i\in C_j} k(Y_i,x_j)\alpha(\mathrm{d}x_j) \tag{28}$$
$$\times \left(\pi \sum_{\boldsymbol{\pi}} \prod_{j=1}^{n(\boldsymbol{\pi})}(n_j-1)! \int_{\mathbb{X}^{n(\boldsymbol{\pi})}} \int_{\mathbb{R}^+} u^{n-1} \mathrm{e}^{-\int_{\mathbb{X}} \log(\beta(s)+u)\alpha_n^{\mathbf{X}}(\mathrm{d}s)}\,\mathrm{d}u \right.$$
$$\left. \times \prod_{i\in C_j} k(Y_i,x_j)\alpha(\mathrm{d}x_j) \right)^{-1}.$$

As a corollary of Proposition 4.2, we obtain the main result of ([33], Theorem 1), namely the posterior density of means of the mixture of Dirichlet process model introduced in [32]. Note that there is a slight inaccuracy in the expression for this density in Corollary 1 of [33], which should coincide with the one below. Moreover, we also obtain a representation for its posterior cumulative distribution function, which is new.

**Corollary 4.1.** *Let $\tilde{Q}$ be the random probability measure associated with the mixture of Dirichlet process. The posterior density function of a mean $\tilde{Q}(g)$, given $\mathbf{Y}$, is then given by*

$$\phi^{\mathbf{Y}}(\sigma;g) = \frac{\sum_{\boldsymbol{\pi}} \prod_{j=1}^{n(\boldsymbol{\pi})}(n_j-1)! \int_{\mathbb{X}^{n(\boldsymbol{\pi})}} \rho^n(\sigma,h) \prod_{i \in C_j} k(Y_i,x_j)\alpha(\mathrm{d}x_j)}{\sum_{\boldsymbol{\pi}} \prod_{j=1}^{n(\boldsymbol{\pi})}(n_j-1)! \int_{\mathbb{X}} \prod_{i \in C_j} k(Y_i,x)\alpha(\mathrm{d}x)}, \qquad (29)$$

*where*

$$\rho^n(\sigma;h) = \frac{a+n-1}{\pi} \int_{\mathbb{R}^+} \mathrm{Re}\{\mathrm{e}^{-\int_{\mathbb{X}} \log(1-\mathrm{i}t(h(s)-\sigma))\alpha_n^{\mathbf{X}}(\mathrm{d}s)}\}\,\mathrm{d}t \qquad (30)$$

*is the posterior density of $\tilde{P}(h)$, with $\tilde{P}$ being the Dirichlet process and, as before, $\alpha_n^{\mathbf{X}}(\cdot) = \alpha(\cdot) + \sum_{j=1}^{n(\boldsymbol{\pi})} n_j \delta_{x_j}$. Moreover, the posterior cumulative distribution function of $\tilde{Q}(g)$, given $\mathbf{Y}$, can be represented as*

$$\mathbb{G}^{\mathbf{Y}}(\sigma;g) = \frac{1}{2} - \frac{\sum_{\boldsymbol{\pi}} \prod_{j=1}^{n(\boldsymbol{\pi})}(n_j-1)! \int_{\mathbb{X}^{n(\boldsymbol{\pi})}} \zeta^n(\sigma;h) \prod_{i \in C_j} k(Y_i,x_j)\alpha(\mathrm{d}x_j)}{\pi \sum_{\boldsymbol{\pi}} \prod_{j=1}^{n(\boldsymbol{\pi})}(n_j-1)! \int_{\mathbb{X}} \prod_{i \in C_j} k(Y_i,x)\alpha(\mathrm{d}x)}, \qquad (31)$$

*with $\zeta^n(\sigma;h) = \int_0^\infty \frac{1}{t} \mathrm{Im}\{\mathrm{e}^{-\int_{\mathbb{X}} \log(1-\mathrm{i}t(h(s)-\sigma))\alpha_n^{\mathbf{X}}(\mathrm{d}s)}\}\,\mathrm{d}t$ being the posterior cumulative distribution function of a Dirichlet mean $\tilde{P}(h)$.*

### 4.2. Generalized gamma NRMI's

Let us now consider the generalized gamma NRMI, which is based on the CRM with intensity

$$\nu(\mathrm{d}v,\mathrm{d}x) = \frac{\gamma}{\Gamma(1-\gamma)} \frac{\mathrm{e}^{-\tau v}}{v^{1+\gamma}}\,\mathrm{d}v\,\alpha(\mathrm{d}x), \qquad (32)$$

where $\gamma \in (0,1)$ and $\tau \geq 0$. This class can be characterized as the exponential family generated by the positive stable laws (see [2, 36]). It includes the stable CRM for $\tau = 0$, the inverse Gaussian CRM for $\gamma = 1/2$ and the gamma CRM as $\gamma \to 0$. Note that the resulting NRMI, termed a *generalized gamma NRMI*, is well defined if and only if $\alpha$ is a finite measure. Before proceeding, recall that $\mathbb{E}[\tilde{P}(\cdot)] = \frac{\alpha(\cdot)}{a} = P_0(\cdot)$, where $a := \alpha(\mathbb{X})$, and $P_0$ is usually referred to as the prior guess at the shape of $\tilde{P}$.

Now, from (6), and by noting that $\psi(g) = \int_{\mathbb{X}}(\tau+g(x))^\gamma \alpha(\mathrm{d}x) - \tau^\gamma a$, it follows immediately that $\tilde{P}(g)$ is finite if and only if

$$\int_{\mathbb{X}}(\tau+t|g(x)|)^\gamma \alpha(\mathrm{d}x) < \infty \qquad \text{for any } t > 0. \qquad (33)$$

Henceforth, we consider functions $g$ such that (33) holds true. In the following proposition, we provide expressions for the prior distribution of $\tilde{P}(g)$.

**Proposition 4.3.** *Let $\tilde{P}$ be a generalized gamma NRMI. The cumulative distribution function can then be expressed as*

$$\mathbb{F}(\sigma;g) = \frac{1}{2} - \frac{e^\beta}{\pi} \int_0^{+\infty} \frac{1}{t} \mathrm{Im}(e^{-\beta \int_{\mathbb{X}}(1-\mathrm{i}t(g(x)-\sigma))^\gamma P_0(\mathrm{d}x)})\,\mathrm{d}t, \qquad (34)$$

*where $\beta = a\tau^\gamma > 0$, or, also, as*

$$\mathbb{F}(\sigma;g) = \frac{1}{2} - \frac{e^\beta}{\pi} \int_0^\infty \frac{1}{t} e^{-\beta A_y(t)} \sin(\beta B_y(t))\,\mathrm{d}t, \qquad (35)$$

*where*

$$A_\sigma(t) = \int_{\mathbb{X}} [1+t^2(g(x)-\sigma)^2]^{\gamma/2} \cos\{\gamma \arctan[t(g(x)-\sigma)]\} P_0(\mathrm{d}x), \qquad (36)$$

$$B_\sigma(t) = \int_{\mathbb{X}} [1+t^2(g(x)-\sigma)^2]^{\gamma/2} \sin\{\gamma \arctan[t(g(x)-\sigma)]\} P_0(\mathrm{d}x). \qquad (37)$$

Let us now turn our attention to posterior quantities. The next result provides both the posterior density and the posterior cumulative distribution function of the mean of a generalized gamma NRMI.

**Proposition 4.4.** *Let $\tilde{P}$ be a generalized gamma NRMI. The posterior density function of $\tilde{P}(g)$, given $\mathbf{X}$, is then of the form (14), with*

$$\chi_g(t,z) = \frac{\gamma \beta^{n(\boldsymbol{\pi})} e^{-\beta \int_{\mathbb{X}}(1-\mathrm{i}t(g(x)-z))^\gamma P_0(\mathrm{d}x)} \prod_{j=1}^{n(\boldsymbol{\pi})}[1-\mathrm{i}t(g(X_j^*)-z)]^{\gamma-n_j}}{\pi \sum_{j=0}^{n-1} \binom{n-1}{j}(-1)^j \beta^{j/\gamma} \Gamma(n(\boldsymbol{\pi})-j/\gamma;\beta)}, \qquad (38)$$

*where $\Gamma(\cdot;\cdot)$ stands for the incomplete gamma function. Moreover, the posterior cumulative distribution function can be written as*

$$\mathbb{F}^{\mathbf{X}}(\sigma;g)$$
$$= \frac{1}{2} - \int_0^\infty \int_0^\infty t^{-1} u^{n-1} \mathrm{Im}\Bigg\{ e^{-\beta \int_{\mathbb{X}}(1+u-\mathrm{i}t(g(x)-\sigma))^\gamma P_0(\mathrm{d}x)}$$
$$\times \prod_{j=1}^{n(\boldsymbol{\pi})}[1+u-\mathrm{i}t(g(X_j^*)-\sigma)]^{\gamma-n_j} \Bigg\} \mathrm{d}u\,\mathrm{d}t$$
$$\times \left( \pi(\gamma \beta^{n(\boldsymbol{\pi})})^{-1} \sum_{j=0}^{n-1} \binom{n-1}{j} (-1)^j \beta^{j/\gamma} \Gamma(n(\boldsymbol{\pi})-j/\gamma,\beta) \right)^{-1}.$$

*Remark 4.1.* At this point, it is useful to compare the expressions arising for extended gamma NRMI's and for generalized gamma NRMI's. As for the latter, by looking at (38) and at the imaginary part of the distribution function of the posterior mean in Proposition 4.4, one can easily identify the characteristic functional of the generalized gamma CRM times a product of $n(\boldsymbol{\pi})$ terms, each of these being the characteristic function of a gamma random variable. Now, each term is clearly associated with a distinct observation $X_j^*$ and the number of times $n_j$ it has been recorded. Moreover, the precise expression arises by taking the derivative of the characteristic functional (which corresponds to observing a certain value for the first time) and then by taking the $(n_j - 1)$th derivative of the resulting expression (which corresponds to observing repetitions of this observation). This structure is apparent in the expressions in Proposition 4.4 and, indeed, it is the idea of seeing the observations as derivatives which inspires the proof of Theorem 3.1. Turning our attention to the extended gamma NRMI which is an extension to the Dirichlet process, one still derives the posterior as derivatives, but, then, due to the important relation

$$\tau_{n_i}(u|y_i)\mathrm{e}^{-\psi(u)} = \frac{\Gamma(n_i)\mathrm{e}^{-\int_{\mathbb{X}} \log(1+u(\beta(x))^{-1})\alpha(\mathrm{d}x)}}{(u+\beta(y_i))^{-n_i}} \qquad (39)$$
$$= \Gamma(n_i)\beta(y_i)^{n_i}\mathrm{e}^{-\int_{\mathbb{X}} \log(1+u(\beta(x))^{-1})\alpha^*(\mathrm{d}x)},$$

one can adjoin the derivatives (observations) to the characteristic functional, getting a sort of quasi-conjugacy, which is a proper conjugacy property if and only if the NRMI is the Dirichlet process. By conjugacy, one usually refers to the fact that the posterior distribution is of the same form as the prior with updated parameters. Here, by quasi-conjugacy, we refer to the fact that, although conjugacy itself does not hold, it does hold conditionally on some latent variable.

We are now in a position to provide the posterior distribution of means of a random mixture density (5) driven by a generalized gamma NRMI. Before stating the result, we introduce the notation $(a)_b = \Gamma(a+b)/\Gamma(a)$ for the Pochhammer symbol.

**Proposition 4.5.** *Let $\tilde{Q}$ be the random probability measure associated with the random density (5) driven by a generalized gamma NRMI. The posterior density function of a mean $\tilde{Q}(g)$, given $\mathbf{Y}$, is then of the form (19), with*

$$\begin{aligned}\xi_h(t,z) = {} & \frac{\gamma}{\pi}\mathrm{e}^{-\beta \int_{\mathbb{X}}(1-\mathrm{i}t(h(x)-z))^{\gamma}P_0(\mathrm{d}x)} \\ & \times \sum_{\boldsymbol{\pi}} \prod_{j=1}^{n(\boldsymbol{\pi})}(1-\gamma)_{n_j-1}\int_{\mathbb{X}}[1-\mathrm{i}t(h(x)-z)]^{\gamma-n_j}\prod_{i\in C_j}k(Y_i,x)\alpha(\mathrm{d}x) \\ & \times \left(\sum_{\boldsymbol{\pi}}\sum_{i=0}^{n-1}\binom{n-1}{i}(-1)^i\beta^{i/\gamma-n(\boldsymbol{\pi})}\Gamma(n(\boldsymbol{\pi})-i/\gamma,\beta)\right. \\ & \left.\qquad \times \prod_{j=1}^{n(\boldsymbol{\pi})}(1-\gamma)_{n_j-1}\int_{\mathbb{X}}\prod_{i\in C_j}k(Y_i,x)\alpha(\mathrm{d}x)\right)^{-1},\end{aligned} \qquad (40)$$

where $h(x) = \int_{\mathbb{Y}} |g(y)| k(y,x) \lambda(\mathrm{d}y)$. Moreover, the posterior distribution function can be represented as (21), with

$$\begin{aligned}
&\zeta_h(\sigma; t) \\
&= \frac{\gamma \beta^{n(\boldsymbol{\pi})}}{\pi} \int_{\mathbb{R}^+} t^{n-1} \mathrm{e}^{-\beta \int_{\mathbb{X}} (1 - \mathrm{i}t(h(x)-z)+u)^\gamma P_0(\mathrm{d}x)} \\
&\quad \times \frac{\sum_{\boldsymbol{\pi}} \prod_{j=1}^{n(\boldsymbol{\pi})} (1-\gamma)_{n_j - 1} \int_{\mathbb{X}} [1 - \mathrm{i}t(h(x) - z) + u]^{\gamma - n_j} \prod_{i \in C_j} k(Y_i, x) \alpha(\mathrm{d}x)\,\mathrm{d}u}{\sum_{\boldsymbol{\pi}} \sum_{i=0}^{n-1} \binom{n-1}{i} (-1)^i \beta^{i/\gamma} \Gamma(n(\boldsymbol{\pi}) - i/\gamma, \beta) \prod_{j=1}^{n(\boldsymbol{\pi})} (1-\gamma)_{n_j - 1} \int_{\mathbb{X}} \prod_{i \in C_j} k(Y_i, x) \alpha(\mathrm{d}x)}.
\end{aligned} \qquad (41)$$

*4.2.1. Connections with the two-parameter Poisson–Dirichlet process*

In this section, we point out a connection between a generalized gamma NRMI and the celebrated two-parameter Poisson–Dirichlet process due to Pitman [35]. The latter random probability measure has found interesting applications in a variety of fields such as population genetics, statistical physics, excursion theory and combinatorics; for details and references, see [37]. A recent systematic study of functionals of the two-parameter Poisson–Dirichlet process is provided in [22]. However, some interesting new results follow from the treatment of the generalized gamma NRMI given above.

Let us first recall the definition of a two-parameter Poisson–Dirichlet random probability measure. Suppose that $\tilde{\mu}_\gamma$ is the $\gamma$-stable CRM arising from (32) when $\tau = 0$. Let $\mathcal{Q}_\gamma$ denote its probability distribution on $(\mathbb{M}, \mathcal{B}(\mathbb{M}))$ and $\mathcal{Q}_{\gamma,\theta}$ be another probability distribution on $(\mathbb{M}, \mathcal{B}(\mathbb{M}))$ such that $\mathcal{Q}_{\gamma,\theta} \ll \mathcal{Q}_\gamma$ and $(\mathrm{d}\mathcal{Q}_{\gamma,\theta}/\mathrm{d}\mathcal{Q}_\gamma)(\mu) = (\mu(\mathbb{X}))^{-\theta}$, where $\theta > -\gamma$. If $\tilde{\mu}_{\gamma,\theta}$ is the random measure whose distribution coincides with $\mathcal{Q}_{\gamma,\theta}$, one then defines the two-parameter Poisson–Dirichlet process as $\tilde{P}_{\gamma,\theta} = \tilde{\mu}_{\gamma,\theta}/\tilde{\mu}_{\gamma,\theta}(\mathbb{X})$. We will now show how to exploit the fact that $\tilde{P}_{\gamma,\theta}(g)$ can be obtained as a suitable mixture of means of generalized gamma NRMI's in order to obtain a new representation for the distribution function of $\tilde{P}_{\gamma,\theta}(g)$.

Denote by $Z$ a gamma random variable with shape parameter $\theta/\gamma > 0$ and scale parameter $\beta > 0$, that is,

$$f_Z(z) = \frac{\beta^{\theta/\gamma}}{\Gamma(\theta/\gamma)} z^{(\theta/\gamma) - 1} \mathrm{e}^{-\beta z}. \qquad (42)$$

Let $\tilde{\mu}_Z$ be a CRM, independent of $Z$, with intensity measure obtained by multiplying the Lévy measure corresponding to the generalized gamma CRM in (32) by $Z$, that is,

$$\nu_Z(\mathrm{d}v, \mathrm{d}x) = \frac{Z\gamma}{\Gamma(1-\gamma)} \frac{\mathrm{e}^{-\tau v}}{v^{1+\gamma}} \mathrm{d}v\, \alpha(\mathrm{d}x),$$

and define $\tilde{P}_Z = \tilde{\mu}_Z/\tilde{\mu}_Z(\mathbb{X})$. By (34), we can set $\tau = 1$ without loss of generality. See also [36], Section 4. It can then be shown that for any function $g: \mathbb{X} \to \mathbb{R}^+$ satisfying the integrability condition (33), one has

$$\int_0^{+\infty} \mathbb{P}[\tilde{P}_z(g) \le x] f_Z(z)\, \mathrm{d}z = \mathbb{P}[\tilde{P}_{\gamma,\theta}(g) \le x]. \qquad (43)$$

The above equality in distribution follows by noting that

$$\mathbb{E}[(\omega + \tilde{P}_Z(g))^{-1}] = \theta \int_0^{+\infty} \frac{\int_{\mathbb{X}}[1+u\omega+ug(x)]^{\gamma-1}\alpha(\mathrm{d}x)}{\{\int_{\mathbb{X}}[1+u\omega+ug(x)]^{\gamma}\alpha(\mathrm{d}x)\}^{(\theta/\alpha)+1}}\,\mathrm{d}u = \mathbb{E}[(\omega+\tilde{P}_{\gamma,\theta}(g))^{-1}]$$

for any $\omega \in \mathbb{C}$ such that $|\arg(\omega)| < \pi$, where the first expected value is computed with respect to the product measure of the vector $(\tilde{P}_Z, Z)$. Using this connection, one can exploit the representation for the distribution of $\tilde{P}_Z(g)$ given in Proposition 4.3 in order to deduce a new, surprisingly simple, expression for the probability distribution of $\tilde{P}_{\gamma,\theta}(g)$, which can be compared with alternative representations given in [22].

**Proposition 4.6.** *Let $g$ be a function for which (33) holds true. The cumulative distribution function of the mean of a two-parameter Poisson–Dirichlet process $\tilde{P}_{\gamma,\theta}(g)$ can then be represented as*

$$\mathbb{F}(\sigma;g) = \frac{1}{2} - \frac{1}{\pi}\int_0^{\infty} \frac{\sin(\theta/\gamma \arctan B_\sigma(t)/A_\sigma(t))}{t[A_\sigma^2(t)+B_\sigma^2(t)]^{\theta/(2\gamma)}}\,\mathrm{d}t, \tag{44}$$

*where $A_\sigma$ and $B_\sigma$ are defined in (36) and (37), respectively.*

As an interesting consequence of the representation in (43), we obtain the finite-dimensional distributions of the two-parameter Poisson–Dirichlet process with $\gamma = 1/2$ and $\theta > 0$, which were first obtained in [3], Theorem 3.1. Before stating the result, it is worth noting that if $\gamma = 1/2$, the generalized gamma CRM reduces to an inverse Gaussian CRM and that, consequently, the finite-dimensional distributions of the two-parameter Poisson–Dirichlet process are obtained as mixtures with respect to those of the inverse Gaussian NRMI.

**Proposition 4.7.** *For any partition of $\mathbb{X}$ into sets $A_1,\ldots,A_n \in \mathscr{X}$ such that $P_0(A_i) = p_i > 0$ for any $i$, a density function of the random vector $(\tilde{P}_{1/2,\theta}(A_1),\ldots,\tilde{P}_{1/2,\theta}(A_{n-1}))$ on the simplex $\Delta_{n-1} = \{(w_1,\ldots,w_{n-1}) \in [0,1]^{n-1} : \sum_{i=1}^{n-1} w_i \leq 1\}$ is given by*

$$f(w_1,\ldots,w_{n-1}) = \frac{(\prod_{i=1}^n p_i)}{\pi^{(n-1)/2}}\frac{\Gamma(\theta+n/2)}{\Gamma(\theta+1/2)}\frac{w_1^{-3/2}\cdots w_{n-1}^{-3/2}(1-\sum_{i=1}^{n-1} w_i)^{-3/2}}{[\mathcal{A}_n(w_1,\ldots,w_{n-1})]^{\theta+n/2}},$$

*where $\mathcal{A}_n(w_1,\ldots,w_{n-1}) = \sum_{i=1}^{n-1} p_i^2 w_i^{-1} + p_n^2(1-\sum_{i=1}^{n-1} w_i)^{-1}$.*

*Remark 4.2.* Two interesting distributional properties of Dirichlet means can be readily extended to means of species sampling models, which include homogeneous NRMI's and the two-parameter Poisson–Dirichlet process as special cases [36]. Recall that a species sampling model is defined as an a.s. discrete random probability measure $\tilde{P}(\cdot) = \sum_{i \geq 1} \tilde{p}_i \delta_{X_i}(\cdot)$ such that the $\tilde{p}_i$'s (weights) are independent from the $X_i$'s (locations), which are i.i.d. from some non-atomic distribution $P_0$. The first property we consider is related to the symmetry of the distribution of $\tilde{P}(g)$. If $\tilde{P}$ is a Dirichlet process, conditions for symmetry have been investigated, for example, in [13, 39]. If $\tilde{P}$ is a species sampling model, it can be shown that, analogously to the Dirichlet case, $\tilde{P}(g)$ is symmetric if the distribution of $g(X)$, where $X$ is a random variable with distribution $P_0$, is symmetric. Another distributional property of a mean of a Dirichlet process was considered in [43], where the author proves that $P_0 \circ g^{-1}$ has the same distribution as $\tilde{P}(g)$ if $P_0 \circ g^{-1}$ is Cauchy distributed. By mimicking the proof of [43], one can easily show that such a property holds true when $\tilde{P}$ is any species sampling model.

## 5. Concluding remarks

In this final section, we briefly discuss two further issues, namely the concrete implementation of the results provided in the paper (also in relation to simulation schemes) and frequentist asymptotics. As for the former, we note that in the numerical computation of the integrals involved in the distributional *formulae* for means of NRMI's, the integrands are typically well behaved, which, in some sense, is quite natural, given the "mean-operation" has a smoothing effect: hence, for specific cases, exact expressions are computable with numerical double-integration packages; see, for example, [39, 40]. For mixtures of NRMI, numerical computation becomes problematic since sums over partitions are involved, which increase the computing time at an exponential rate. As an alternative, one could also resort to simulation schemes to obtain an approximation of the distribution of interest. However, the simulation of realizations of an NRMI is a delicate task, since it is based on CRM's which jump infinitely often on any bounded interval. Consequently, a simulation algorithm is necessarily based on some truncation, which, in some cases, may compromise posterior inferences; see, for example, [10]. This implies that the availability of exactly computable expressions for distributions of means is also useful for practical purposes as a benchmark for testing the accuracy of a simulation algorithm: one simulates trajectories of the process (by tuning the number of jumps and the number of trajectories) until a suitable distance between the exact and the simulated distribution of the mean is less than a prescribed error. If this is the case, the numerical output can then be exploited to compute any quantity of interest (not just means).

In order to evaluate frequentist asymptotic properties of Bayesian procedures, the data are assumed to be independently generated by a "true" distribution $P_{\mathrm{tr}}$. As far as consistency is concerned, one essentially has that NRMI's are consistent if $P_{\mathrm{tr}}$ is discrete and inconsistent if $P_{\mathrm{tr}}$ is non-atomic (with the exception of the Dirichlet process). This can be informally illustrated by looking at, for instance, the $\gamma$-stable NRMI: the predictive distributions associated with such an NRMI are of the form

$$\mathbb{P}(X_{n+1} \in \cdot | X_1^*, \ldots, X_{n(\boldsymbol{\pi})}^*) = \frac{\gamma n(\boldsymbol{\pi})}{n} P_0(\cdot) + \frac{1}{n} \sum_{i=1}^{n(\boldsymbol{\pi})} (n_i - \gamma) \delta_{X_i^*}(\cdot). \qquad (45)$$

If the data are generated independently from a non-atomic $P_{\mathrm{tr}}$, all observations are distinct, which implies that $n(\boldsymbol{\pi}) = n$ almost surely. Hence, as the sample size diverges, the predictive distribution converges weakly, almost surely, to $\gamma P_0 + (1-\gamma) P_{\mathrm{tr}}$ and $\mathbb{E}[\tilde{P}(g)|\mathbf{X}]$ converges to $\gamma P_0(g) + (1-\gamma) P_{\mathrm{tr}}(g)$. Consequently, the posterior distribution of $\tilde{P}(g)$ is

inconsistent unless $P_0 = P_{\text{tr}}$. On the other hand, if $P_{\text{tr}}$ is discrete, then $\lim n(\boldsymbol{\pi})/n = 0$ almost surely and $\mathbb{E}[\tilde{P}(g)|\mathbf{X}]$ converges to the "true" value $P_{\text{tr}}(g)$; see [21, 24] for general results in this direction. Inconsistency with respect to a continuous $P_{\text{tr}}$ is not really problematic since, for modelling continuous data, one would naturally use an NRMI mixture rather than an NRMI itself and, in such a case, we have consistency under mild assumptions [30]. A related issue of interest is the validity of Bernstein–von Mises theorems, which provide information about asymptotic normality of the posterior distribution and asymptotic equivalence with respect to MLE estimation. For the nonparametric case, only a few results, both positive and negative, are currently known; see [16] for a stimulating account of the topic. The derivation of such results can be a formidable task since the posterior representation of the random probability measure is explicitly involved. This implies that no result is known, even for the mixture of Dirichlet process. Recently, in [21], a Bernstein–von Mises-type theorem has been obtained for the two-parameter Poisson–Dirichlet process in the case of continuous $P_{\text{tr}}$. Under the same assumptions as in Theorem 3.1 of [21], we can then deduce the asymptotic behavior for the posterior distribution of a mean functional of the two-parameter Poisson–Dirichlet process, denoted by $\tilde{P}_{\gamma,\theta}^{\mathbf{X}}(g)$. Indeed, we obtain that as the number of observations generated by an absolutely continuous $P_{\text{tr}}$ diverges,

$$n^{1/2}\{\tilde{P}_{\gamma,\theta}^{\mathbf{X}}(g) - \mathbb{E}[\tilde{P}_{\gamma,\theta}^{\mathbf{X}}(g)]\}$$
$$\xrightarrow{d} \mathrm{N}(0, (1-\gamma)[P_{\text{tr}}(g^2) - [P_{\text{tr}}(g)]^2] + \gamma(1-\gamma)[P_0(g^2) - [P_0(g)]^2]$$
$$+ \gamma(P_{\text{tr}}(g) - P_0(g))^2) \qquad \text{a.s.,}$$

where $\mathrm{N}(\mu, \sigma^2)$ stands for a Gaussian random variable with mean $\mu$ and variance $\sigma^2$. In particular, if $g = \mathbb{1}_A$, then the asymptotic behavior of $\tilde{P}_{\gamma,\theta}^{\mathbf{X}}(A)$ is obtained, that is, the centered and rescaled distribution of $\tilde{P}_{\gamma,\theta}^{\mathbf{X}}(A)$ converges to a mean zero Gaussian random variable with variance $(1-\gamma)[p(1-p)] + \gamma(1-\gamma)[q(1-q)] + \gamma(p-q)^2$, having set $p := P_{\text{tr}}(A)$ and $q := P_0(A)$. In order to derive such results for general NRMI's, an approach similar to [21] seems advisable, although the technicalities may become overwhelming. This may well represent a topic for future research.

## Appendix

**Proof of Theorem 3.1.** In order to derive a representation of the posterior density, we start by discretizing $\tilde{P}$ according to the procedure of Regazzini and Sazonov [41]. It essentially consists of discretizing the random probability measure and the sample along a tree of nested partitions of $\mathbb{X}$ which, at level $m$, is made up of sets $B_{m,1},\ldots,B_{m,k_m+1}$ with $B_{m,k_m+1} = \mathbb{X} \setminus (\bigcup_{i=1}^{k_m} B_{m,i})$ such that $B_{m,k_m+1} \downarrow \varnothing$ and $\max_{1 \le i \le k_m} \text{diam}(B_{m,i}) \to 0$ as $m$ tends $+\infty$, where $\text{diam}(B)$ is the diameter of $B$. The discretized random mean, at level $m$ of the tree, will be of the form $\tilde{P}_m(g) = \sum_{j=1}^{k_m+1} g(b_{m,j})\tilde{P}(B_{m,j})$, where $b_{m,j}$ is any point in $B_{m,j}$ for $j = 1,\ldots,k_m+1$. Moreover, denote by $\mathscr{M}_m$ the marginal distribution of $\mathbf{X}$ with respect to $\tilde{P}_m$. Whenever the $j$th distinct element, $x_j^*$, lies in $B_{m,i}$, it is as if we had observed $b_{m,i}$; see [41] for details. Note that whatever tree of partitions has been chosen, there always exists an $m^*$ such that for every $m > m^*$, the $n(\boldsymbol{\pi})$ distinct observations within the sample fall in $n(\boldsymbol{\pi})$ distinct sets of the partition. Given such a discretization, Propositions 3 and 4 in [40] easily extend to the case of NRMI's over Polish spaces. Thus, the discretized posterior mean converges, in the sense of almost sure weak convergence, to the actual posterior mean. If the observations $\mathbf{X}$ are such that $n_{s_j}$ of them lie in $B_{m,s_j}$, for $j = 1, \ldots, n(\boldsymbol{\pi})$, an expression of the posterior density function of the discretized mean $\tilde{P}_m(g)$ is given by

$$\frac{(-1)^{n+1}}{\mathscr{M}_m(\mathbf{X})} \frac{\partial^n}{\partial r_{m,s_1}^{n_{s_1}} \cdots \partial r_{m,s_{n(\boldsymbol{\pi})}}^{n_{s_{n(\boldsymbol{\pi})}}}} \tag{46}$$
$$\times I_{a^+}^{n-1} \mathbb{F}(\sigma; r_{m,0}, \ldots, r_{m,k_m+1}) \bigg|_{(r_{m,1},\ldots,r_{m,k_m+1})=(f(b_{m,1}),\ldots,f(b_{m,k_m+1}))},$$

where $\mathbb{F}$ stands for the prior distribution of the discretized mean and $I_{a^+}^n h(\sigma) = \int_a^\sigma \frac{(\sigma-u)^{n-1}}{(n-1)!} h(u)\, du$ is the *Liouville–Weyl fractional integral*, for $n \geq 1$, whereas $I_{a^+}^0$ represents the identity operator. We now move on to the identification of the limiting density. To this end, set $\alpha_{m,j} = \alpha(B_{m,j})$, for $j = 1, \ldots, k_{m+1}$, and rewrite (46) as

$$\frac{(-1)^{n+1}}{\pi \mathscr{M}_m(\mathbf{X})} I_{a^+}^{n-1} \operatorname{Im} \int_0^{+\infty} \frac{1}{t} e^{-\sum_{j=0}^{k_{m+1}} \int_0^{+\infty}(1-e^{itv(r_{m,j}-\sigma)})\rho(dv|b_{m,j})\alpha_{m,j}} \prod_{l=1}^{n(\boldsymbol{\pi})} \Lambda_{\alpha_{m,s_l}}^{n_{s_l}}(t)\, dt, \tag{47}$$

where $\Lambda$ is defined as

$$\Lambda_{\alpha_{m,s_l}}^{n_{s_l}}(t) := e^{\int_0^{+\infty}(1-e^{itv(r_{m,j}-\sigma)})\rho(dv|b_{m,s_l})\alpha_{m,s_l}} \left\{ \frac{\partial^{n_{s_l}}}{\partial r_{m,s_l}^{n_{s_l}}} e^{-\int_0^{+\infty}(1-e^{itv(r_{m,j}-\sigma)})\rho(dv|b_{m,s_l})\alpha_{m,s_l}} \right\}.$$

By virtue of the diffuseness of $\alpha$, one has

$$\Lambda_{\alpha_{m,s_l}}^{n_{s_l}}(t) = (it)^{n_{s_l}} \alpha_{m,s_l} \int_0^{+\infty} v^{n_{s_l}} e^{itv(r_{m,s_l}-\sigma)} \rho(dv|b_{m,s_l}) + o(\alpha_{m,s_l})$$

as $m \to \infty$. In order to complete the proof, one needs to resort to the expression of the marginal distribution of the observations provided in Proposition 4 of [23]. Hence, if we let $m$ tend to $+\infty$ and apply Theorem 35.7 in Billingsley [1] and dominated convergence, the desired result follows. Note that, as a by-product, we have also proven that the posterior distribution of the means is absolutely continuous with respect to the Lebesgue measure on $\mathbb{R}$.

The proof of the representation of the posterior cumulative distribution function in (16) consists of the following steps. First, we use the idea suggested by [13], that is,

$$\mathbb{F}^{\mathbf{X}}(\sigma; g) = \mathbb{P}\{\tilde{\mu}(g - \sigma \mathbf{1}) \leq 0 | \mathbf{X}\}. \tag{48}$$

As in [40], we now resort to Gurland's inversion formula, which, combined with Proposition 2 in [40], yields

$$\mathbb{F}^{\mathbf{X}}(\sigma;g) = \frac{1}{2} - \frac{1}{\pi}\lim_{T\uparrow+\infty}\int_0^T \frac{1}{t}\operatorname{Im}\{\mathbb{E}[\mathrm{e}^{\mathrm{i}t\tilde{\mu}(g-\sigma\mathbf{1})}|\mathbf{X}]\}\,\mathrm{d}t.$$

The conclusion can now be easily deduced from [23], Theorem 1, according to which

$$\mathbb{E}[\mathrm{e}^{\mathrm{i}t\tilde{\mu}(g-\sigma\mathbf{1})}|\mathbf{X}] = \int_0^{+\infty}\mathbb{E}[\mathrm{e}^{\mathrm{i}t\tilde{\mu}^{(u)}(g-\sigma\mathbf{1})+\mathrm{i}t\sum_{r=1}^{n(\boldsymbol{\pi})}(g(x_r^*)-\sigma)J_r^{(u,\mathbf{X})}}]f_{U_n}^{\mathbf{X}}(u)\,\mathrm{d}u. \qquad \square$$

**Proof of Theorem 3.2.** The first thing to note is that by (7), determining the distribution of $(\tilde{Q}(g)|\mathbf{Y})$ is equivalent to determining the distribution of $(\tilde{P}(h)|\mathbf{Y})$ with $h(x) = \int_{\mathbb{Y}} g(y)k(y,x)\lambda(\mathrm{d}y)$. Moreover, by linearity of the mean, one has $(\tilde{P}(h)|\mathbf{Y}) = \int_{\mathbb{X}} h(x) \times (\tilde{P}(\mathrm{d}x)|\mathbf{Y})$. Thus, we need a posterior representation of the NRMI, given the data $\mathbf{Y}$ which come from the mixture of NRMI in (5). To this end, one can adapt Theorem 2 in [18] to obtain the disintegration

$$(\tilde{P}|\mathbf{Y}) \stackrel{d}{=} \int_{\mathbb{X}^n}(\tilde{P}|\mathbf{X})\mathscr{M}(\mathrm{d}\mathbf{X}|\mathbf{Y}), \qquad (49)$$

where $\mathscr{M}(\mathbf{X}|\mathbf{Y})$ stands for the marginal distribution of the latent variables $\mathbf{X}$, given the observables $\mathbf{Y}$. From (49) combined with Theorem 3.1, a first description of the distribution of $(\tilde{Q}(g)|\mathbf{Y})$ follows. Indeed, the posterior density can be represented as

$$\phi^{\mathbf{Y}}(\sigma;g) = \int_{\mathbb{X}^n}\rho^{\mathbf{X}}(\sigma;h)\mathscr{M}(\mathrm{d}\mathbf{X}|\mathbf{Y}), \qquad (50)$$

where $\rho^{\mathbf{X}}(\sigma;h)$ coincides with the density function given in Theorem 3.1 with $h(x) = \int_{\mathbb{Y}} g(y)k(y,x)\lambda(\mathrm{d}y)$. Moreover, the posterior distribution function is of the form

$$\mathbb{G}^{\mathbf{Y}}(\sigma;g) = \int_{\mathbb{X}^n}\mathbb{F}^{\mathbf{X}}(\sigma;h)\mathscr{M}(\mathrm{d}\mathbf{X}|\mathbf{Y}), \qquad (51)$$

where $\mathbb{F}^{\mathbf{X}}(\sigma;h)$ is given in Theorem 3.1, with $h$ defined as above. Note that the previous expressions could also have been obtained by combining the discretization of the observation space employed in [33] with the limiting arguments used in the proof of Theorem 3.1.

In order to get explicit descriptions of the posterior distributions in (50) and (51), we need to find expressions for the marginal distribution of $\mathbf{X}$ given $\mathbf{Y}$. From Bayes' theorem, it follows that

$$\mathscr{M}(\mathrm{d}\mathbf{X}|\mathbf{Y}) = \frac{\prod_{i=1}^n k(Y_i,X_i)\mathscr{M}(\mathrm{d}\mathbf{X})}{\int_{\mathbb{X}^n}\prod_{i=1}^n k(Y_i,X_i)\mathscr{M}(\mathrm{d}\mathbf{X})}. \qquad (52)$$

It is clear that the marginal distribution of $\mathbf{X}$ can be described by the joint distribution of the distinct variables $\mathbf{X}^*$ and the induced random partition $\boldsymbol{\pi}$. In Proposition 4 of [23], an expression for the joint distribution of $(\mathbf{X}^*, \boldsymbol{\pi})$ is provided. By inserting this formula into (52), suitably rearranging the terms and making use of the partition notation introduced at the beginning of Section 2, one obtains

$$\mathscr{M}(\mathrm{d}\mathbf{X}|\mathbf{Y}) = \frac{\prod_{j=1}^{n(\boldsymbol{\pi})} \prod_{i \in C_j} k(Y_i, X_j^*) \alpha(\mathrm{d}X_j^*) [\int_{\mathbb{R}^+} u^{n-1} \mathrm{e}^{-\psi(u)} \prod_{j=1}^{n(\boldsymbol{\pi})} \tau_{n_j}(u|X_j^*) \, \mathrm{d}u]}{\int_{\mathbb{R}^+} u^{n-1} \mathrm{e}^{-\psi(u)} [\sum_{\boldsymbol{\pi}} \prod_{j=1}^{n(\boldsymbol{\pi})} \int_{\mathbb{X}} \prod_{i \in C_j} k(Y_i, x) \tau_{n_j}(u|x) \alpha(\mathrm{d}x)] \, \mathrm{d}u}, \quad (53)$$

where $\sum_{\boldsymbol{\pi}}$ stands for the sum over partitions and the $\tau_{n_j}$'s are defined in (11). Inserting (53) into (50) and carrying out suitable simplifications, one obtains the desired expression for the posterior density. As for the posterior distribution function, insert (53) into (51) and, after some algebra, the result follows. □

**Proof of Proposition 4.1.** The proofs of both representations are based on Theorem 3.1 and on the key relation (39) discussed in Remark 4.1 of Section 4.2. As for the posterior density function, note that $\tau_{n_j}(u|X_j^*) = \Gamma(n_j)[\beta(X_j^*) + u]^{-n_j}$ and that $\kappa_{n_j}(\mathrm{i}t[g(X_j^*) - z]|X_j^*) = \Gamma(n_j)[\beta(X_j^*) - \mathrm{i}t(g(X_j^*) - z)]^{-n_j}$ for $j = 1, \ldots, n(\boldsymbol{\pi})$. Inserting these expression into the general one of Theorem 3.1 and carrying out the appropriate simplifications using (39) allows the posterior density function of $\tilde{P}(g)$ to be written as desired. With reference to the posterior cumulative distribution function, one has that, given $\mathbf{X}$ and $U_n$, the $j$th jump $J_j^{(U_n, \mathbf{X})}$ is gamma distributed with parameters $(\beta(X_j^*) + U_n, n_j)$, for $j = 1, \ldots, n(\boldsymbol{\pi})$. Hence, we have $\zeta_g(\sigma; u, t) = \exp[-\int_{\mathbb{X}} \log(1 - \mathrm{i}t(g(x) - \sigma)[\beta(x) + u]^{-1}) \alpha^*(\mathrm{d}x)]$ for (17). Finally, it is easy to verify that $t^{-1}\zeta_g(\sigma; u, t)$ is absolutely integrable in $(M, +\infty)$ for any $M > 0$ and $u > 0$. Thus, appropriate simplifications between denominator and numerator lead to the result. □

**Details for the determination of (24).** Given $g = \mathbb{1}_A$ and $\alpha(A) = \alpha(A^c) = 1$, the integral in (23) can be written as

$$\int_0^\infty \frac{\cos(2 \arctan t(1-z)/\beta_1 - \arctan tz/\beta_2)}{[\beta_1^2 + t^2(1-z)^2][\beta_2^2 + t^2z^2]^{1/2}} \, \mathrm{d}t$$

$$= \int_0^\infty \frac{1}{[\beta_1^2 + t^2(1-z)^2][\beta_2^2 + t^2z^2]^{1/2}}$$

$$\times \left\{ \cos^2\left(\arctan \frac{t(1-z)}{\beta_1}\right) \cos\left(\arctan \frac{tz}{\beta_2}\right) - \sin^2\left(\arctan \frac{t(1-z)}{\beta_1}\right) \cos\left(\arctan \frac{tz}{\beta_2}\right) \right.$$

$$\left. + 2\sin\left(\arctan \frac{t(1-z)}{\beta_1}\right) \cos\left(\arctan \frac{t(1-z)}{\beta_1}\right) \sin\left(\arctan \frac{tz}{\beta_2}\right) \right\} \mathrm{d}t$$

and it can be easily seen that it reduces to

$$\int_0^\infty \frac{1}{[\beta_1^2 + t^2(1-z)^2][\beta_2^2 + t^2z^2]^{1/2}}$$

$$\times \left\{ \frac{\beta_1^2}{\beta_1^2 + t^2(1-z)^2} \frac{\beta_2}{\sqrt{\beta_2^2 + t^2 z^2}} \right.$$

$$\left. - \frac{t^2(1-z)^2}{\beta_1^2 + t^2(1-z)^2} \frac{\beta_2}{\sqrt{\beta_2^2 + t^2 z^2}} + 2 \frac{t(1-z)\beta_1}{\beta_1^2 + t^2(1-z)^2} \frac{tz}{\sqrt{\beta_2^2 + t^2 z^2}} \right\} dt$$

$$= \int_0^\infty \frac{t^2(1-z)(2\beta_1 z - \beta_2(1-z)) + \beta_1^2 \beta_2}{[\beta_1^2 + t^2(1-z)^2]^2 [\beta_2^2 + t^2 z^2]} \, dt$$

$$= (1-z)(2\beta_1 z - \beta_2(1-z)) \frac{1}{(1-z)^4 z^2} \frac{\pi}{4\beta_1/(1-z)(\beta_1/(1-z) + \beta_2/z)^2}$$

$$+ \frac{\pi \beta_1^2 \beta_2}{2(1-z)^4 z^2} \left\{ \frac{z}{\beta_2 (\beta_1^2/(1-z)^2 - \beta_2^2/z^2)^2} \right.$$

$$\left. - \frac{1}{2(\beta_1^2/(1-z)^2 - \beta_2^2/z^2)\beta_1^3/(1-z)^3} - \frac{1-z}{(\beta_1^2/(1-z)^2 - \beta_2^2/z^2)^2 \beta_1} \right\}$$

$$= \frac{\pi[2\beta_1 z + \beta_2(1-z)]}{4\beta_1[\beta_1 z + \beta_2(1-z)]^2} + \frac{\pi[2\beta_1^3 z^3 - 3\beta_1^2 \beta_2 z^2(1-z) + \beta_2^3(1-z)^3]}{4\beta_1[\beta_1^2 z^2 - \beta_2^2(1-z)^2]^2}$$

$$= \frac{\pi z}{[\beta_1 z + \beta_2(1-z)]^2},$$

from which (24) follows. $\square$

**Proof of Proposition 4.2.** The representations in (25) and (27) are obtained by a direct application of Theorem 3.2: one only needs to compute the quantities already exploited in the proof of Proposition 4.1 and insert them into the relevant expressions of Theorem 3.2. In order to obtain the posterior representations in (26) and (28) which make use of the quasi-conjugacy of the extended gamma process, it is enough to apply Fubini's theorem and use (39). $\square$

**Proof of Corollary 4.1.** Recall that the Dirichlet case corresponds to the extended gamma case with $\beta(x) = c > 0$, which, without loss of generality, we can set equal to 1. We first derive the posterior density function. To this end, consider (26) and note that in this case, within the denominator, one has

$$\int_{\mathbb{R}^+} u^{n-1} e^{-\log(1+u)(a+n)} \, du = \frac{\Gamma(a)\Gamma(n)}{\Gamma(a+n)}, \tag{54}$$

having set $\alpha(\mathbb{X}) := a$. This allows (26) to be rewritten as

$$\xi_h(t, z) = \sum_{\boldsymbol{\pi}} \int_{\mathbb{X}^{n(\boldsymbol{\pi})}} \frac{\Gamma(a+n) e^{-\int_\mathbb{X} \log(\beta(s) - it(h(s) - z))\alpha_n^\mathbf{X}(ds)}}{\pi \Gamma(a)\Gamma(n)}$$

$$\times \frac{\prod_{j=1}^{n(\boldsymbol{\pi})} (n_j - 1)! \prod_{i \in C_j} k(Y_i, x_j) \alpha(dx_j)}{\sum_{\boldsymbol{\pi}} \prod_{j=1}^{n(\boldsymbol{\pi})} (n_j - 1)! \int_\mathbb{X} \prod_{i \in C_j} k(Y_i, x) \alpha(dx)}, \tag{55}$$

which, combined with (20) and an application of Fubini's theorem, leads to a posterior density of the form

$$\phi^{\mathbf{Y}}(\sigma;g) = \frac{\sum_{\boldsymbol{\pi}} \prod_{j=1}^{n(\boldsymbol{\pi})}(n_j-1)! \int_{\mathbb{X}^{n(\boldsymbol{\pi})}} \varrho^n(\sigma,h) \prod_{i\in C_j} k(Y_i,x_j)\alpha(\mathrm{d}x_j)}{\sum_{\boldsymbol{\pi}} \prod_{j=1}^{n(\boldsymbol{\pi})}(n_j-1)! \int_{\mathbb{X}} \prod_{i\in C_j} k(Y_i,x)\alpha(\mathrm{d}x)}, \qquad (56)$$

where $\varrho^n$ is the posterior density of $\tilde{P}(h)$, given $\mathbf{X}$, as given in [23], Theorem 1. But, since the posterior mean of a Dirichlet process is again a Dirichlet process mean with updated parameter measure $\alpha_n^{\mathbf{X}}$, one can choose the simplest possible expression for such a density, which, to date, is (30) obtained in [40]. Hence, $\varrho^n$ can be replaced by (30) in (56), and (29) follows. Note that (30) can also be derived from the posterior density in Proposition 4.1 by considering the case $n=1$, for which (54) yields $a^{-1}$. To obtain (31), one proceeds in a similar fashion, starting from (28), and then simplifies the resulting expression using the representation for the cumulative distribution function of Dirichlet process mean provided in [39]. □

**Proof of Proposition 4.3.** We start by noting that

$$\mathbb{P}(\tilde{P}(g) \leq \sigma) = \mathbb{P}(\tilde{\mu}(\tau[g-\sigma\mathbf{1}]) \leq 0) \qquad (57)$$

and then apply Gurland's inversion formula to the right-hand side of (57). Showing that $\mathbb{F}$ is absolutely continuous (with respect to the Lebesgue measure on $\mathbb{R}$), one has that the prior distribution of $\tilde{P}(g)$ is of the form (34), with $\beta = a\tau^\gamma > 0$. Note that, thanks to the reparameterization induced by (57), the generalized gamma NRMI is completely specified by the parameters $P_0$, $\gamma$ and $\beta$. In order to make (34) completely explicit, introduce the quantities (36) and (37) and observe that

$$\exp\left\{-\beta \int_{\mathbb{X}} [1-\mathrm{i}t(g(x)-\sigma)]^\gamma P_0(\mathrm{d}x)\right\} = \exp\{-[\beta A_\sigma(t) - \mathrm{i}\beta B_\sigma(t)]\}.$$

This implies that (34) can be represented, by working out the imaginary part, as (35). □

**Proof of Proposition 4.4.** In order to obtain the expressions for the posterior distribution, we resort to Theorem 3.1, but we have to take into account the reparametrization in (57): hence, we set, without loss of generality, $\tau = 1$ and $\alpha(\mathrm{d}x) = \beta P_0(\mathrm{d}x)$ in (32). Consequently, we get $\tau_{n_j}(u|X_j^*) = \gamma[1+u]^{\gamma-n_j}(1-\gamma)_{n_j-1}$ and $\kappa_{n_j}(\mathrm{i}t(g(X_j^*)-z)|X_j^*) = \gamma[1-\mathrm{i}t(g(X_j^*)-z)]^{\gamma-n_j}(1-\gamma)_{n_j-1}$, for $j=1,\ldots,n(\boldsymbol{\pi})$, where, as before, $(a)_b$ stands for the Pochhammer symbol. The explicit form of the denominator is obtained by a suitable change of variable combined with the binomial theorem. Some algebra leads to an analytic representation of the posterior density of $\tilde{P}(g)$ as in (14) with the function $\mathcal{X}_g$ as in (38). This proves that first part of the proposition. In order to derive the posterior cumulative distribution function given in (16), first note that the jumps $J_i^{(U_n,\mathbf{X})}$ are gamma distributed with scale parameter $U_n + \gamma$ and shape parameter $n_i - \gamma$. Then, write the explicit form of the other quantities involved and carry out some simplifications. It is then easy to see that $\zeta(\sigma; u, t)$ is absolutely integrable in $(M, +\infty)$ for any $M > 0$ and $u > 0$, thus, the result follows. □

**Proof of Proposition 4.5.** The result is obtained by exploiting the quantities computed for deriving Proposition 4.4 and inserting them into the expression of Theorem 3.2. This, combined with some algebra, leads to the formulae in (40) and (41). □

**Proof of Proposition 4.6.** Let $\mathbb{F}^Z(\cdot; g)$ denote the cumulative distribution function of the mean $\tilde{P}_Z(g)$, given $Z$. By Proposition 4.3, $\mathbb{F}^Z(\cdot; g)$ coincides with (34), with $z\beta$ in place of $\beta$. Reasoning as in the proof of Proposition 4.3, we obtain

$$\mathbb{F}^Z(\sigma; g) = \frac{1}{2} - \frac{e^{z\beta}}{\pi} \int_0^\infty \frac{1}{t} e^{-z\beta A_\sigma(t)} \sin(z\beta B_\sigma(t)) \, dt,$$

where $A_\sigma$ and $B_\sigma$ are defined as in (36) and (37), respectively. We now integrate with respect to $Z$ in order to obtain the unconditional cumulative distribution function of the two-parameter Poisson–Dirichlet process, that is,

$$\mathbb{F}(\sigma; g) = \frac{1}{2} - \frac{\beta^{\theta/\gamma}}{\pi \Gamma(\theta/\gamma)} \int_0^\infty dt \frac{1}{t} \int_0^\infty z^{\theta/\gamma - 1} e^{-z\beta A_\sigma(t)} \sin(z\beta B_\sigma(t)) \, dz$$

$$= \frac{1}{2} - \frac{\beta^{\theta/\gamma}}{\pi} \int_0^\infty \frac{\sin(\theta/\gamma \arctan B_\sigma(t)/A_\sigma(t))}{t[\beta^2 A_\sigma^2(t) + \beta^2 B_\sigma^2(t)]^{\theta/(2\gamma)}} \, dt,$$

where the last equality follows from 3.944.5 in [15]. Simplifying with respect to $\beta$ leads to the desired result. □

**Proof or Proposition 4.7.** From Proposition 1 in [29], the density function of the random vector $(\tilde{P}_Z(A_1), \ldots, \tilde{P}_Z(A_{n-1}))$ with $\gamma = 1/2$ is given by

$$f^Z(w_1, \ldots, w_{n-1}) = \frac{e^{Z\beta}(Z\beta)^n \prod_{i=1}^n p_i}{2^{n/2-1} \pi^{n/2}} K_{-n/2}(Z\beta \sqrt{\mathcal{A}_n(w_1, \ldots, w_{n-1})})$$

$$\times w_1^{-3/2} \cdots w_{n-1}^{-3/2} \left(1 - \sum_{i=1}^{n-1} w_i\right)^{-3/2} \{(Z\beta)^2 \mathcal{A}_n(w_1, \ldots, w_{n-1})\}^{-n/4}.$$

If one integrates the above density with respect to $f_Z$ in (42) and makes use of 6.561.16 in [15], the result easily follows. □

## Acknowledgements

The authors are grateful to two anonymous referees for their valuable comments and suggestions. L.F. James was supported in part by Grants RGC-HKUST 600907, SBM06/07.BM14 and HIA05/06.BM03 of the HKSAR. A. Lijoi and I. Prünster were partially supported by MIUR Grant 2008MK3AFZ.